\documentclass[a4paper]{article}

\usepackage{amssymb}
\usepackage{amsmath}
\usepackage{amsthm}

\usepackage{endnotes}

\newcommand{\ssection}[1] { 
\refstepcounter{section}
\setcounter{equation}{0} 
\setcounter{subsection}{0} 
\addcontentsline{toc}{section}{\textbf{\thesection.\ #1}} 
\bigskip\bigskip\noindent
\Large
\textbf{\thesection.\ #1}
\nopagebreak\bigskip\nopagebreak
\normalsize
} 
\def\thesection{{\arabic{section}}}

\newcounter{appendice}

\newtheorem{theorem}{Theorem}[section]
\newtheorem{lemma}[theorem]{Lemma}

\newtheorem{corollary}[theorem]{Corollary}

\theoremstyle{definition}
\newtheorem{definition}[theorem]{Definition}

\theoremstyle{definition}

\theoremstyle{remark}
\newtheorem{remark}[theorem]{Remark}

\newcommand{\G}{\Gamma}
\newcommand{\B}{\text{B}}
\newcommand{\tg}{\tilde{\gamma}^{(A)}_q}
\newcommand{\tb}{\tilde{\beta}^{(A)}_q}

\newcommand{\R}{\mathbb R}	\newcommand{\Z}{\mathbb Z}
\newcommand{\C}{\mathbb C}
 
\begin{document}


\par\vskip 1cm\vskip 2em

\begin{center}

{\Huge On integral representations of $q$-gamma and $q$--beta functions}

\par

\vskip 2.5em \lineskip .5em
 
{\large
\begin{tabular}[t]{c}
$\mbox{Alberto De Sole}\!\!\phantom{m}\mbox{,\ \ Victor G. Kac}$\\
\end{tabular}
\par
}
\medskip
{\small
\begin{tabular}[t]{ll}
& Department of Mathematics, MIT \\
& 77 Massachusetts Avenue, Cambridge, MA 02139, USA \\
& E--mails: {\tt desole@math.mit.edu} \\
& \hspace{1.25cm} {\tt kac@math.mit.edu} 
\end{tabular}
}
\bigskip
\end{center}


\vskip 1 em
\centerline{\bf Abstract} 
\smallskip
{\small 
\noindent
We study $q$--integral representations of the $q$--gamma and the $q$--beta 
functions. This study leads to a very interesting $q$--constant.
As an application of these integral representations,
we obtain a simple conceptual proof of 
a family of identities for Jacobi triple product, including Jacobi's identity,
and of Ramanujan's formula for the bilateral hypergeometric series.
}

\bigskip


\renewcommand{\theenumi}{\alph{enumi}}
\renewcommand{\labelenumi}{(\theenumi)}

\ssection{Introduction}

The $q$--gamma function $\G_q(t)$, a $q$--analogue of Euler's gamma function,
was introduced by Thomae \cite{thomae} and later by Jackson 
\cite{jackson_1904} as the infinite product
\begin{equation}\label{int1}
\G_q(t)=\frac{(1-q)_q^{t-1}}{(1-q)^{t-1}}\ ,\quad t>0\ ,
\end{equation}
where $q$ is a fixed real number $0<q<1$. Here and further we use the
following notation:
\begin{eqnarray}
(a+b)_q^n &=& \prod_{j=0}^{n-1}(a+q^jb)\ ,\quad \text{ if } n\in\Z_+\ ,
 \label{intr2}\\
(1+a)_q^\infty &=& \prod_{j=0}^{\infty}(1+q^ja)\ ,
 \label{intr3}\\
(1+a)_q^t &=& \frac{(1+a)_q^\infty}{(1+q^ta)_q^\infty}\ ,
 \quad \text{ if } t\in\C\ . \label{intr4}
\end{eqnarray}
Notice that, under our assumptions on $q$, the infinite product (\ref{intr3})
is convergent. Moreover, the definitions (\ref{intr2}) and (\ref{intr4}) 
are consistent.

Though the literature on the $q$--gamma function and its applications 
is rather extensive, \cite{exton}, \cite{gasper_rahman}, \cite{aar},
the authors usually avoided the use
of its $q$--integral representation.
In fact, each time when a $q$--integral representation was discussed, it was,
as a rule, not quite right.
The first correct integral representation of $\G_q(t)$ that we know of
is in reference \cite{koornwinder_92}:
\begin{equation}\label{intr5}
\G_q(t)=\int_0^{\frac{1}{1-q}}x^{t-1}E_q^{-qx}d_qx\ .
\end{equation}
Here $E_q^x$ is one of the two $q$--analogues of the exponential function:
\begin{eqnarray}
E_q^x &=& \sum_{n=0}^\infty q^{n(n-1)/2}\frac{x^n}{[n]!}\ 
=\ (1+(1-q)x)_q^\infty\ , \label{intr6} \\
e_q^x &=& \sum_{n=0}^\infty \frac{x^n}{[n]!}\ 
=\ \frac{1}{(1-(1-q)x)_q^\infty}\ , \label{intr7}
\end{eqnarray}
and the $q$--integral (introduced by Thomae \cite{thomae} 
and Jackson \cite{jackson_1910}) is defined by
\begin{equation}\label{intr8}
\int_0^{a}f(x)d_qx
=(1-q)\sum_{j=0}^\infty aq^jf(aq^j)\ .
\end{equation}

The $q$--beta function was more fortunate in this respect.
Already in the mentioned papers by Thomae and Jackson it was shown
that the $q$--beta function defined by the usual formula
\begin{equation}\label{intr9}
\B_q(t,s)=\frac{\G_q(s)\G_q(t)}{\G_q(s+t)}\ ,
\end{equation}
has the $q$--integral representation, which is a $q$--analogue 
of Euler's formula:
\begin{equation}\label{intr10}
\B_q(t,s)=\int_0^1x^{t-1}(1-qx)_q^{s-1}d_qx\ ,\quad t,s>0\ .
\end{equation}

Jackson \cite{jackson_1910} made an attempt to give a $q$--analogue
of another Euler's integral representation of the beta function:
\begin{equation}\label{intr11}
\B(t,s)
 = \int_0^\infty \frac{x^{t-1}}{(1+x)^{t+s}}dx\ .
\end{equation}
However, his definition is not quite right, since it is not quite equal
to $B_q(t,s)$, as will be explained in Remark \ref{errore}.
A correct $q$--analogue of (\ref{intr11}) is the famous Ramanujan's formula
for the bilateral hypergeometric series, see \cite[pp 502--505]{aar}
(in fact, Ramanujan's formula was known already to Kronecker).

In the present paper we give a $q$--integral representation of $\G_q(t)$
based on the $q$--exponential function $e_q^x$,
and give a $q$--integral representation of $\B_q(t,s)$ which is a 
$q$--analogue of (\ref{intr11}).
Both representations are based on the following remarkable function:
\begin{equation}\label{intr12}
K(x,t)
=\frac{x^t}{1+x}\Big(1+\frac{1}{x}\Big)_q^t(1+x)_q^{1-t}\ .
\end{equation}
This function is a $q$--constant in $x$, i.e.
$$
K(qx,t)=K(x,t)\ ,
$$
and for $t$ an integer it is indeed independent on $x$, and is equal
to $q^{t(t-1)/2}$.
However, for $t\in(0,1)$ this function does depend on $x$, since for these $t$
one has
$$
\lim_{q\rightarrow0}K(x;t)=x^t+x^{t-1}\ .
$$
Our integral representations are as follows:
\begin{eqnarray}
\G_q(t) 
&=& K(A,t)\int_0^{\infty/A(1-q)}x^{t-1}e_q^{-x}d_qx\ , \label{intr13} \\
\B_q(t,s)
&=& K(A,t)\int_0^{\infty/A}\frac{x^{t-1}}{(1+x)_q^{t+s}}d_qx\ , \label{intr14}
\end{eqnarray}
where the improper integral, following \cite{jackson_1910}
and \cite{koornwinder_99}, is defined by
\begin{equation}\label{intr15}
\int_0^{\infty/A}f(x)d_qx
=(1-q)\sum_{n\in\Z}\frac{q^n}{A}f\Big(\frac{q^n}{A}\Big)\ .
\end{equation}
Since $K(A,t)$ depends on $A$, we conclude that integrals in both formulas
do depend on $A$.
In his formulas Jackson used the factor $q^{t(t-1)/2}$ in place
of $K(A,t)$, which is correct only for an integer $t$.

Note also that formula (\ref{intr5}) for $\G_q(t)$ can be written via an
improper integral too
(since $E_q^{-\frac{q^{n}}{1-q}}=0$ for $n\leq0$):
\begin{equation}\label{intr16}
\G_q(t)=\int_0^{\infty/(1-q)}x^{t-1}E_q^{-qx}d_qx\ .
\end{equation}

We will refer to the book \cite{kac_qc}
for notation and basic facts on $q$--calculus.
Unfortunately the authors of the book didn't know about the reference
\cite{koornwinder_92} and made the same mistake 
as their predecessors in the definition (\ref{intr16}) (taking $\infty$
in place of $\infty/(1-q)$, which gives a divergent series).
However all their arguments hold verbatim for the definition
(\ref{intr5}) (or (\ref{intr16})) and are used in the present paper
to derive (\ref{intr13}) and (\ref{intr14}).

In Sections 5 and 6 we will apply equation (\ref{intr14})
to find an integral representation of the $q$--beta function
which is manifestly symmetric under the exchange of $t$ and $s$,
and to find a $q$--analogue of translation invariance of certain improper
integrals.
Finally, in Section 7 we will show that 
equation (\ref{intr13}) is equivalent to a family of triple product
identities, a special case of which is the Jacobi triple product identity:
\begin{equation}\label{intr17a}
(1-q)_q^\infty(1-x)_q^\infty(1-q/x)_q^\infty
=\sum_{n=-\infty}^\infty(-1)^nq^{n(n-1)/2}x^n\ ,
\end{equation}
equation (\ref{intr14}) is equivalent to the famous Ramanujan's identity,
see \cite[pp 501--505]{aar}:
\begin{equation}\label{intr17b}
\sum_{n=-\infty}^\infty\frac{(1-a)_q^n}{(1-b)_q^n}x^n
=\frac{(1-q)_q^\infty(1-b/a)_q^\infty(1-ax)_q^\infty(1-q/ax)_q^\infty}{(1-b)_q^\infty(1-q/a)_q^\infty(1-x)_q^\infty(1-b/ax)_q^\infty}\ ,
\end{equation}
and the symmetric integral representation of the $q$--beta function
is equivalent to the following identity:
\begin{equation}\label{intr17c}
\sum_{n=-\infty}^\infty
\frac{(1-a)_q^n(1-q/a)_q^{-n}}{(1-b)_q^n(1-c)_q^{-n}}
=\frac{(1-a)_q^\infty(1-q/a)_q^\infty}{(1-b)_q^\infty(1-c)_q^\infty}
\frac{(1-q)_q^\infty(1-bc/q)_q^\infty}{(1-b/a)_q^\infty(1-ac/q)_q^\infty}\ .
\end{equation}

One of the authors wishes to thank A. Varchenko for sending corrections
to the book \cite{kac_qc}, including the one mentioned above.


\ssection{Notation and preliminary results}


Throughout this paper we will assume $q$ to be a fixed number between 0 and 1.
We denote by $D_q$ the $q$--derivative of a function
$$
(D_q f)(x)=\frac{f(qx)-f(x)}{(q-1)x}\ .
$$
The Jackson definite integral of the function $f$ is defined to be 
\cite{jackson_1910}
$$
\int_0^a f(x)d_qx=(1-q)\sum_{n\geq0}q^naf(q^na)\ .
$$
Notice that the series on the right--hand side is guaranteed to be convergent
as soon as the function $f$ is such that, for some $C>0,\ \alpha>-1$,
$|f(x)|<Cx^{\alpha}$ in a right neighborhood of $x=0$.
Jackson integral and $q$--derivative are related by the ``fundamental
theorem of quantum calculus'', \cite[p 73]{kac_qc}
\begin{theorem}
\begin{enumerate}
\item If $F$ is any anti $q$--derivative of the function $f$, 
namely $D_qF=f$, continuous at $x=0$, then
$$
\int_0^af(x)d_qx=F(a)-F(0)\ .
$$
\item For any function $f$ one has
$$
D_q\int_0^xf(t)d_qt=f(x)\ .
$$
\end{enumerate}
\end{theorem}
\begin{remark}
It is easy to check the $q$--analogue of the Leibniz rule:
$$
D_q(f(x)g(x))=g(x)D_qf(x)+f(qx)D_qg(x)\ .
$$
An immediate consequence is the $q$--analogue 
of the rule of integration by parts:
$$
\int_0^ag(x)D_qf(x)d_qx
=f(x)g(x)\Big|_0^a-\int_0^af(qx)D_qg(x)d_qx\ .
$$
\end{remark}


One defines the Jackson integral in a generic interval $[a,b]$ by 
\cite{jackson_1910}:
$$
\int_a^bf(x)d_qx=\int_0^bf(x)d_qx-\int_0^af(x)d_qx\ .
$$
One also defines improper integrals in the following way
\cite{jackson_1910}, \cite{koornwinder_99}:
\begin{equation}\label{improp}
\int_0^{\infty/A}f(x)d_qx
=(1-q)\sum_{n\in\Z}\frac{q^n}{A}f\Big(\frac{q^n}{A}\Big)\ .
\end{equation}
\begin{remark}
Notice that in order the series on the right--hand side to be convergent,
it suffices that the function $f$ satisfies the conditions:
$|f(x)|<Cx^\alpha, \quad \forall x\in [0,\epsilon)$, 
for some $C>0,\ \alpha>-1,\ \epsilon>0$,
and $|f(x)|<Dx^\beta,\quad\forall x\in [N,\infty)$, 
for some $D>0,\ \beta<-1,\ N>0$.
In general though, even when these conditions are satisfied,
the value of the sum will be dependent on the constant $A$.
In order the integral to be independent of $A$, the anti $q$--derivative
of $f$ needs to have limits for $x\rightarrow 0$ and $x\rightarrow+\infty$.
\end{remark}


One has the following {\itshape reciprocity relations}:
\begin{eqnarray}\label{reciprocity}
&& \int_0^Af(x)d_qx
=\int_{q/A}^{\infty/A}\frac{1}{x^2}f\Big(\frac{1}{x}\Big)d_qx\ , \nonumber\\
&& \int_0^{\infty/A}f(x)d_qx
=\int_0^{\infty\cdot A}\frac{1}{x^2}f\Big(\frac{1}{x}\Big)d_qx\ . 
\end{eqnarray}
This is a special case of the following more general change 
of variable formula, \cite[p 107]{kac_qc}.
If $u(x)=\alpha x^\beta$, then
$$
\int_{u(a)}^{u(b)}f(u)d_qu
=\int_a^bf(u(x))D_{q^{1/\beta}}u(x)d_{q^{1/\beta}}x\ .
$$


The $q$--analogue of an integer number $n$ is
$$
[n]=\frac{1-q^n}{1-q}
=1+q+\dots+q^{n-1}\ .
$$
In general we will denote $[t]=\frac{1-q^t}{1-q}$ even for a non integer $t$.


\renewcommand{\theenumi}{\arabic{enumi}}
\renewcommand{\labelenumi}{(\theenumi)}

The following lemma follows immediately from the definitions
(\ref{intr2}), (\ref{intr3}) and (\ref{intr4}) of $(a+b)_q^t$.
The proof is left as an exercise to the reader. 
See also \cite[pp 106-107]{kac_qc}.
\begin{lemma}\label{prelim}
Let $n\in\Z_+$ and $t,s,a,b,A,B\in\R$.
\begin{enumerate}
\item $D_qx^t=[t]x^{t-1}$
\item $D_q(Ax+b)_q^n=[n]A(Ax+b)_q^{n-1}$
\item $D_q(a+Bx)_q^n=[n]B(a+Bqx)_q^{n-1}$
\item $D_q(1+Bx)_q^t=[t]B(1+Bqx)_q^{t-1}$
\item $D_q\frac{Ax^s}{(1+Bx)_q^t}
=[s]\frac{Ax^{s-1}}{(1+Bx)_q^{t+1}}-B([t]-[s])\frac{Ax^s}{(1+Bx)_q^{t+1}}$
\item $D_q\frac{(1+Ax)_q^s}{(1+Bx)_q^t}
=[s]A\frac{(1+Aqx)_q^{s-1}}{(1+Bqx)_q^t}-B[t]\frac{(1+Ax)_q^s}{(1+Bx)_q^{t+1}}$
\item $(1+x)_q^{s+t}=(1+x)_q^s(1+q^sx)_q^t$
\item $(1+x)_q^{-t}=\frac{1}{(1+q^{-t}x)_q^t}$
\item $(1+q^sx)_q^t
=\frac{(1+x)_q^{s+t}}{(1+x)_q^s}
=\frac{(1+q^tx)_q^s}{(1+x)_q^s}(1+x)_q^t$
\item $(1+q^{-n}x)_q^t
=\frac{(x+q)_q^n}{(q^tx+q)_q^n}(1+x)_q^t$
\end{enumerate}
\end{lemma}

\renewcommand{\theenumi}{\alph{enumi}}
\renewcommand{\labelenumi}{(\theenumi)}

The $q$--analogues of the exponential function are given by (\ref{intr6})
and (\ref{intr7}).
The equality between the series expansion and the infinite product
expansion of $e_q^x$ and $E_q^x$ (in the domain where both expansions converge)
is proved by taking the limit for $n\rightarrow\infty$
in the Gauss and Heine's $q$--binomial formulas, \cite[pp 29--32]{kac_qc}.

\begin{remark}
Notice that for $q\in(0,1)$ the series expansion of $e_q^x$ has 
radius of convergence $1/(1-q)$.
This corresponds to the fact that the infinite product 
$1/(1-(1-q)x)_q^\infty$ has a pole at $x=1/(1-q)$.
On the contrary, the series expansion of $E_q^x$ converges for every $x$.
Both product expansions (\ref{intr6}) and (\ref{intr7})  converge for all $x$.
\end{remark}

\begin{lemma}\label{exp_prop}
The $q$--exponential functions satisfy the following properties:
\begin{enumerate}
\item $D_qe_q^x=e_q^x\ ,\qquad D_qE_q^x=E_q^{qx}$.
\item $e_q^xE_q^{-x}=E_q^xe_q^{-x}=1$.
\end{enumerate}
\end{lemma}
\noindent The proof is straightforward and it is left as an exercise 
to the reader. See also \cite[pp 29--32]{kac_qc}.


\ssection{Definition of $q$--gamma and $q$--beta functions}

The Euler's gamma and beta functions are defined as the following
definite integrals ($s,t>0$):
\begin{eqnarray}
\G(t) &=& \int_0^\infty x^{t-1}e^{-x}dx\ , \label{cl_gamma} \\
\B(t,s) &=& \int_0^1 x^{t-1}(1-x)^{s-1}dx\ , \label{cl_beta1} \\
 &=& \int_0^\infty \frac{x^{t-1}}{(1+x)^{t+s}}dx\ . \label{cl_beta2}
\end{eqnarray}
Notice that (\ref{cl_beta2}) follows immediately from (\ref{cl_beta1})
after a change of variable $x=1/(1+y)$.
From the expression (\ref{cl_beta1}) it is clear that $\B(t,s)$ is symmetric
in $t$ and $s$.
Recall some of the main properties of the gamma and beta functions:
\begin{eqnarray}
&& \G(t+1)=t\G(t)\ ,\quad \G(1)=1\ , \label{cl_prop1}\\
&& \B(t,s)=\G(t)\G(s)/\G(t+s)\ \label{cl_prop2}.
\end{eqnarray}

In this paper we are interested in the $q$--analogue of the gamma and beta
functions.
They are defined in the following way. 
\begin{definition}
\begin{enumerate}
\item For $t>0$, the $q$--gamma function is defined to be
\begin{equation}\label{big_gamma}
\G_q(t)=\int_0^{1/(1-q)}x^{t-1}E_q^{-qx}d_qx\ .
\end{equation}
\item For $s,t>0$, the $q$--beta function is
\begin{equation}\label{big_beta}
\B_q(t,s)=\int_0^1x^{t-1}(1-qx)_q^{s-1}d_qx\ .
\end{equation}
\end{enumerate}
\end{definition}


$\G_q(t)$ and $\B_q(t,s)$ are the ``correct'' $q$--analogues of the
gamma and beta functions, since they reduce to $\G(t)$ and $\B(t,s)$ 
respectively in the limit $q\rightarrow1$, and they satisfy properties
analogues to (\ref{cl_prop1}) and (\ref{cl_prop2}).
This is stated in the following
\begin{theorem}\label{big}
\begin{enumerate}
\item $\G_q(t)$ can be equivalently expressed as
\begin{equation}\label{a}
\G_q(t)=\frac{(1-q)_q^{t-1}}{(1-q)^{t-1}}\ .
\end{equation}
In particular one has
$$
\G_q(t+1)=[t]\G_q(t)\ ,\quad \forall t>0\ ,\quad \G_q(1)=1\ .
$$
\item The $q$--gamma and $q$-beta functions are related to each other 
by the following two equations
\begin{eqnarray}
\G_q(t) &=& \frac{\B_q(t,\infty)}{(1-q)^t}\ , \label{b1} \\
\B_q(t,s) &=& \frac{\G_q(t)\G_q(s)}{\G_q(t+s)}\ . \label{b2}
\end{eqnarray}
\end{enumerate}
\end{theorem}
\begin{proof}
We reproduce here the proof of Kac and Cheung \cite[pp 76--79]{kac_qc},
because similar arguments will be used to
prove the results in the next section.
If we put $s=\infty$ in the definition of the $q$--beta
function, use (\ref{intr6})
and the change of variable $x=(1-q)y$, we get
\begin{eqnarray*}
\B_q(t,\infty) &=& \int_0^1x^{t-1}E_q^{-\frac{qx}{1-q}}d_qx \\
 &=& \frac{1}{(1-q)^t}\int_0^{1/(1-q)}y^{t-1}E_q^{-qy}d_qy\ 
=\ \frac{1}{(1-q)^t}\G_q(t)\ ,
\end{eqnarray*}
which proves (\ref{b1}).
It follows from $q$--integration by parts and
Lemma \ref{prelim} (parts (3) and (4))
that $\B_q(t,s)$ satisfies the following recurrence relations ($t,s>0$):
\begin{eqnarray*}
\B_q(t+1,s) &=& \frac{[t]}{[s]}\B_q(t,s+1)\ , \\
\B_q(t,s+1) &=& \B_q(t,s)-q^s\B_q(t+1,s)\ .
\end{eqnarray*}
Putting these two conditions together we get
$$
\B_q(t,s+1)=\frac{[s]}{[s+t]}\B_q(t,s)\ .
$$
Since clearly $\B_q(t,1)=\frac{1}{[t]}$, we get,
for $t>0$ and any positive integer $n$,
\begin{eqnarray}
\B_q(t,n) &=& \frac{[n-1]\dots[1]}{[t+n-1]\dots[t]} \nonumber
 \ =\  (1-q)\frac{(1-q)_q^{n-1}}{(1-q^t)_q^n} \\
 &=& (1-q)\frac{(1-q)_q^{n-1}(1-q)_q^{t-1}}{(1-q)_q^{t+n-1}}\ .\label{1}
\end{eqnarray}
Taking the limit for $n\rightarrow\infty$ in this expression we get
$$
\B_q(t,\infty)=(1-q)(1-q)_q^{t-1}\ .
$$
This together with (\ref{b1}) proves (\ref{a}).
We are left to prove (\ref{b2}).
By comparing (\ref{1}) and (\ref{a}) we have that (\ref{b2}) is true
for any positive integer value of $s$.
To conclude that (\ref{b2}) holds for non integer values of $s$
we will use the following simple argument.
If we substitute $a=q^s$ and $b=q^t$ in (\ref{b2}) we can write 
the left--hand side as
$$
(1-q)\sum_{n\geq0}b^n\frac{(1-q^n)_q^\infty}{(1-aq^{n-1})_q^\infty}\ ,
$$
and the right--hand side as
$$
(1-q)\frac{(1-q)_q^\infty(1-ab)_q^\infty}{(1-a)_q^\infty(1-b)_q^\infty}\ .
$$
Both these expressions can be viewed as formal power series in $q$
with coefficients rational functions in $a$ and $b$.
Since we already know that they coincide, for any given $b$,
for infinitely many values of $a$ (of the form $a=q^n$, with positive 
integer $n$), it follows that they must be equal for every value 
of $a$ and $b$.
This concludes the proof of the Theorem.
\end{proof}


\ssection{An equivalent definition of $q$--gamma and $q$--beta functions}


In the previous section the definition of $\G_q(t)$ was obtained from the
integral expression (\ref{cl_gamma}) of the Euler's gamma function, 
simply by replacing the integral with a Jackson integral and 
the exponential function $e^{-x}$ with its $q$--analogue $E_q^{-qx}$.
It is natural to ask what happens if we use the other $q$--exponential
function. In other words, we want to study the following function ($A>0$):
\begin{equation}\label{little_gamma}
\tg(t)=\int_0^{\infty/A(1-q)}x^{t-1}e_q^{-x}d_qx\ .
\end{equation}
Similarly, the function $\B_q(t,s)$ was obtained by taking the $q$--analogue
of the integral expression (\ref{cl_beta1}) of the Euler's beta function.
We now want to study the $q$--analogue of the integral expression
(\ref{cl_beta2}). We thus define
\begin{equation}\label{little_beta}
\tb(t,s)=\int_0^{\infty/A}\frac{x^{t-1}}{(1+x)_q^{t+s}}d_qx\ .
\end{equation}
In this section we will show how the functions $\tg(t)$ and $\tb(t,s)$
are related to the $q$--gamma and $q$--beta function respectively.


We want to adapt the arguments in the proof of Theorem \ref{big} to the
functions $\tg(t)$ and $\tb(t,s)$.
First, by taking the limit $s\rightarrow\infty$ in the definition
of $\tb(t,s)$, using the infinite product expansion of $e_q^x$
and making the change of variables $x=(1-q)y$, we get
\begin{eqnarray*}
\tb(t,\infty) &=& \int_0^{\infty/A}\frac{x^{t-1}}{(1+x)_q^\infty}d_qx\
=\ \int_0^{\infty/A}x^{t-1}e_q^{-\frac{x}{1-q}} d_qx\\
 &=& (1-q)^t\int_0^{\infty/A(1-q)}y^{t-1}e_q^{-y}d_qy\
=\ (1-q)^t\tg(t)\ .
\end{eqnarray*}
We therefore proved
\begin{equation}\label{little_b1}
\tg(t)=\frac{1}{(1-q)^t}\tb(t,\infty)\ .
\end{equation}

We now want to find recursive relations for $\tg(t)$ and $\tb(t,s)$.
By integration by parts we get
$$
\tg(t+1)
=q^{-t}[t]\tg(t)\ .
$$
Here we used the fact that $x^te_q^{-x}$ tends to zero as $x\rightarrow0$
and $x\rightarrow+\infty$ (The second fact follows 
from Lemma \ref{exp_prop} (b)).
Since obviously $\tg(1)=1$, we conclude that for every positive integer $n$
(and any value of $A>0$),
\begin{equation}\label{gamma_integer}
q^{n(n-1)/2}\tg(n)=[n-1]!=\G_q(n)\ .
\end{equation}
Let us now consider the function $\tb(t,s)$.
From integration by parts and the results in Lemma \ref{prelim}
we get ($t,s>0$)
\begin{eqnarray}\label{5}
\tb(t+1,s) 
 &=& -\frac{1}{[t+s]}q^{-t}\int_0^{\infty/A}(qx)^t
    D_q\frac{1}{(1+x)_q^{t+s}}d_qx \\
 &=& q^{-t}\frac{1}{[t+s]}\int_0^{\infty/A}
    \frac{1}{(1+x)_q^{t+s}}D_qx^td_qx\ =\ q^{-t}\frac{[t]}{[t+s]}\tb(t,s)\ .
\nonumber
\end{eqnarray}
For $t=1$ we have
\begin{equation}\label{6}
\tb(1,s)=\int_0^{\infty/A}\frac{1}{(1+x)_q^{s+1}}d_qx=\frac{1}{[s]}\ .
\end{equation}
Formulas (\ref{5}) and (\ref{6}) imply ($s>0,\ n\in\Z_+$)
\begin{equation}\label{beta_int_s}
q^{n(n-1)/2}\tb(n,s)
=(1-q)\frac{(1-q)_q^{s-1}(1-q)_q^{n-1}}{(1-q)_q^{s+n-1}}=\B_q(n,s)\ .
\end{equation}
Similarly we have
\begin{eqnarray}\label{8}
\tb(t,s+1) 
 &=& \frac{1}{[t+s]}q^{s}\int_0^{\infty/A}\frac{1}{(qx)^s}
    D_q\frac{x^{t+s}}{(1+x)_q^{t+s}}d_qx \\
 &=& -q^s\frac{1}{[t+s]}\int_0^{\infty/A}
    \frac{x^{t+s}}{(1+x)_q^{t+s}}D_q\frac{1}{x^s}d_qx\ 
=\ \frac{[s]}{[t+s]}\tb(t,s)\ .\nonumber
\end{eqnarray}
We now need to compute $\tb(t,1)$. By definition and Lemma \ref{prelim} (5)
\begin{equation}\label{9}
\tb(t,1)=\int_0^{\infty/A}\frac{x^{t-1}}{(1+x)_q^{t+1}}d_qx
=\frac{1}{[t]}\int_0^{\infty/A}D_q\frac{x^t}{(1+x)_q^t}d_qx\ .
\end{equation}
When using the fundamental theorem of $q$--calculus to compute 
the right--hand side, we have to be careful, since the limit 
for $x\rightarrow +\infty$ of the function $F(x)=\frac{x^t}{(1+x)_q^t}$ 
does not exist.
On the other hand, by definition of $q$--derivative and Jackson integral,
we have
$$
\int_0^{\infty/A}D_qF(x)d_qx
=\lim_{N\rightarrow\infty}F\Big(\frac{1}{Aq^N}\Big)
-\lim_{N\rightarrow\infty}F\Big(\frac{q^N}{A}\Big)\ ,
$$
where the limits on the right--hand side are taken over the sequence 
of integer numbers $N$.
We then have from (\ref{9})
\begin{equation}\label{4.6}
\tb(t,1)=\frac{1}{[t]}
\Big(\lim_{N\rightarrow\infty}(Aq^N)^t(1+\frac{1}{Aq^N})_q^t\Big)^{-1}\ .
\end{equation}
If we denote by $K(A;t)$ the limit in parenthesis in the right--hand side
of (\ref{4.6}), we can use Lemma \ref{prelim} (10) to get
\begin{eqnarray*}
K(A;t) &=& A^t
  \lim_{N\rightarrow\infty}q^{Nt}\Big(1+\frac{q^{-N}}{A}\Big)_q^t \\
 &=& A^t\Big(1+\frac{1}{A}\Big)_q^t
  \lim_{N\rightarrow\infty}q^{Nt}
  \frac{\Big(\frac{1}{A}+q\Big)_q^N}{\Big(\frac{q^t}{A}+q\Big)_q^N} \\
 &=& A^t\Big(1+\frac{1}{A}\Big)_q^t
  \lim_{N\rightarrow\infty}
  \frac{(1+qA)_q^N}{(1+q^{1-t}A)_q^N} \\
 &=& \frac{1}{1+A}A^t\Big(1+\frac{1}{A}\Big)_q^t(1+A)_q^{1-t}\ .
\end{eqnarray*}
From (\ref{8}) and (\ref{4.6}) we conclude that 
for any $t>0$ and positive integer $n$
\begin{equation}\label{beta_int_t}
K(A;t)\tb(t,n)
=(1-q)\frac{(1-q)_q^{n-1}(1-q)_q^{t-1}}{(1-q)_q^{n+t-1}}=\B_q(t,n)\ .
\end{equation}

In the following lemma we enumerate some interesting properties of 
the function 
$$
K(x;t)=\frac{1}{1+x}x^t\Big(1+\frac{1}{x}\Big)_q^t(1+x)_q^{1-t}\ .
$$
\begin{lemma}\label{prop_k}
\begin{enumerate}
\item In the limit $q\rightarrow1$ and $0$ we have
\begin{eqnarray*}
&& \lim_{q\rightarrow1}K(x;t)=1\ ,\qquad \forall\ x,t\in\R \\
&& \lim_{q\rightarrow0}K(x;t)=x^t+x^{t-1}\ ,
 \qquad \forall\ t\in(0,1),\ x\in\R\ .
\end{eqnarray*}
In particular $K(x,t)$ is not constant in $x$.
\item Viewed as a function of $t$, $K(x;t)$ satisfies the following
recurrence relation:
$$
K(x;t+1)=q^tK(x;t)\ .
$$
Since obviously $K(x;0)=K(x;1)=1$, we have in particular that
for any positive integer $n$
$$
K(x;n)=q^{n(n-1)/2}\ .
$$
\item As function of $x$, $K(x;t)$ is a ``$q$--constant'', namely
$$
D_qK(x,t)=0\ ,\quad \forall\ t,x\in\R\ .
$$
In other words $K(q^nx;t)=K(x;t)$ for every integer $n$.
\end{enumerate}
\end{lemma}
\begin{proof}
The limit for $q\rightarrow1$ of $K(x;t)$ is obviously 1.
In the limit $q\rightarrow0$ we have, for any $\alpha>0$,
$$
(1+x)_q^\alpha=(1+x)\frac{(1+qx)_q^\infty}{(1+q^\alpha x)_q^\infty}\
\longrightarrow\ (1+x)\ .
$$
We therefore have, for $t\in(0,1)$:
$\lim_{q\rightarrow0}K(x;t)
=x^t\Big(1+\frac{1}{x}\Big)$.
For part (b), it follows from the definition of $K(x;t)$ 
and Lemma \ref{prelim} that 
\begin{eqnarray*}
K(x;t+1) &=& \frac{1}{1+x}x^{t+1}\Big(1+\frac{1}{x}\Big)_q^{t+1}(1+x)_q^{-t} \\
 &=& x\Big(1+\frac{q^t}{x}\Big)\frac{1}{(1+q^{-t}x)}K(x;t)\
 =\ q^tK(x;t)\ .
\end{eqnarray*}
For part (c) it suffices to prove that $K(qx;t)=K(x;t)$.
By definition
$$
K(qx;t) = \frac{1}{1+qx}(qx)^t\Big(1+\frac{1}{qx}\Big)_q^t(1+qx)_q^{1-t}\ .
$$
We can replace in the right--hand side
\begin{eqnarray*}
&\Big(1+\frac{1}{qx}\Big)_q^t
=\frac{1+\frac{1}{qx}}{1+\frac{q^t}{qx}}\Big(1+\frac{1}{x}\Big)_q^t&\ , \\
&(1+qx)_q^{1-t}
=\frac{1+q^{1-t}x}{1+x}(1+x)_q^{1-t}&\ .
\end{eqnarray*}
The claim follows from the following trivial identity
$$
\frac{q^t\Big(1+\frac{1}{qx}\Big)(1+q^{1-t}x)}{(1+qx)\Big(1+\frac{1}{q^{1-t}x}\Big)}=1\ .
$$
This concludes the proof of the lemma.
\end{proof}
\begin{remark}
The function $K(x;t)$ is an interesting example of a function which 
is not constant in $x$ 
and with $q$--derivative identically zero.
\end{remark}

It follows from (\ref{beta_int_s}), (\ref{beta_int_t}) and Lemma \ref{prop_k}
that the functions $K(A;t)\tb(t,s)$ and $\B_q(t,s)$
coincide for any $A>0$ as soon as either $t$ or $s$ is a positive integer.
We want to prove that they actually coincide for every $t,s>0$.
\begin{theorem}\label{main}
For every $A,t,s>0$ one has:
\begin{eqnarray}
K(A;t)\tg(t) &=& \G_q(t)\ , \label{t1}\\
K(A;t)\tb(t,s) &=& \B_q(t,s) \label{t2}\ .
\end{eqnarray}
\end{theorem}
\begin{remark}\label{errore}
This result corrects and generalizes a similar statement 
of Jackson \cite{jackson_1910}.
There (\ref{t2}) is proved in the special case in which $s+t$ 
is a positive integer, But, due to a computational mistake, the factor $K(A;t)$
is missing.
\end{remark}
\begin{proof}
(\ref{t1}) is an immediate corollary of (\ref{b1}), (\ref{little_b1})
and (\ref{beta_int_t}).
As in the proof of Theorem \ref{big}, in order to prove (\ref{t2})
it suffices to prove that $K(A;t)\tb(t,s)$ can be written as formal power
series in $q$ with coefficients rational functions in $a=q^s$ and $b=q^t$.
After performing a change of variable $y=Ax$, we get
\begin{equation}\label{2}
K(A;t)\tb(t,s)
=\frac{1}{1+A}\Big(1+\frac{1}{A}\Big)_q^t(1+A)_q^{1-t}
\int_0^{\infty/1}\frac{y^{t-1}}{\Big(1+\frac{y}{A}\Big)_q^{t+s}}d_qy\ .
\end{equation}
Fix $A>0$. After letting $b=q^t$, we can rewrite the factor in front 
of the integral as
$$
\frac{1}{1+A}
\frac{\Big(1+\frac{1}{A}\Big)_q^\infty}{\Big(1+\frac{b}{A}\Big)_q^\infty}
\frac{(1+A)_q^\infty}{\Big(1+\frac{qA}{b}\Big)_q^\infty}\ ,
$$
and this is manifestly a formal power series in $q$ with coefficients rational
functions in $b$.
We then only need to study the integral term in (\ref{2}),
which we decompose as
\begin{equation}\label{decomp}
\int_0^1\frac{y^{t-1}}{\Big(1+\frac{y}{A}\Big)_q^{t+s}}d_qy
+\int_1^{\infty/1}\frac{y^{t-1}}{\Big(1+\frac{y}{A}\Big)_q^{t+s}}d_qy
\end{equation}
After letting $a=q^s$ and $b=q^t$ the first term in (\ref{decomp})
can be written as
$$
(1-q)\sum_{n\geq0}
b^n\frac{\Big(1+ab\frac{q^n}{A}\Big)}{\Big(1+\frac{q^n}{A}\Big)}\ ,
$$
and this is a formal power series in $q$
with coefficients rational functions in $a$ and $b$.
Consider now the second term of (\ref{decomp}).
By relation (\ref{reciprocity}) we can rewrite it as
\begin{equation}\label{bad}
\int_0^q\frac{x^{s-1}}{x^{t+s}\Big(1+\frac{1}{Ax}\Big)_q^{t+s}}d_qx\ .
\end{equation}
Recalling the definition of $K(x;t)$, we have the identity
$$
\frac{1}{x^{t+s}\Big(1+\frac{1}{Ax}\Big)_q^{t+s}}
=\frac{1}{1+Ax}(1+Ax)_q^{1-t-s}\frac{A^{t+s}}{K(Ax;t+s)}\ .
$$
The main observation is that, even though $K(Ax;t+s)$ is not constant 
in $x$, by Lemma \ref{prop_k} $K(Aq^n;t+s)=K(A;t+s),\ \forall\ n\in\Z$,
therefore inside the Jackson integral it can be treated as a constant.
Using this fact, we can rewrite (\ref{bad}) as
\begin{equation}\label{final}
\frac{A^{t+s}}{K(A;t+s)}\int_0^q
\frac{1}{1+Ax}x^{s-1}(1+Ax)_q^{1-t-s}d_qx\ .
\end{equation}
After letting $a=q^s$ and $b=q^t$ we can finally rewrite the first factor
in (\ref{final}) as
\begin{equation}\label{f1}
\frac{(1+A)\Big(1+\frac{ab}{A}\Big)_q^\infty\Big(1+\frac{qA}{ab}\Big)_q^\infty}
{\Big(1+\frac{1}{A}\Big)_q^\infty\Big(1+A\Big)_q^\infty}
\end{equation}
and the integral term in (\ref{final}) as
\begin{equation}\label{f2}
(1-q)\sum_{n\geq0}a^{n+1}
\frac{(1+Aq^{n+2})_q^\infty}
{\Big(1+\frac{Aq^{n+2}}{ab}\Big)_q^\infty}\ .
\end{equation}
Clearly both expression (\ref{f1}) and (\ref{f2}) are
formal power series in $q$ with coefficients 
rational functions in $a$ and $b$.
This concludes the proof of the theorem.
\end{proof}


\ssection{Application 1: an integral expression of the $q$--beta function
manifestly symmetric in $t$ and $s$}

Theorem $\ref{big}$ implies that $\B_q(t,s)$ is a symmetric function
in $t$ and $s$. This is not obvious from 
its integral expression (\ref{big_beta}).
We now want to use Theorem \ref{main} to find an integral expression 
for $\B_q(t,s)$ which is manifestly symmetric under the exchange 
of $t$ and $s$.
By Theorem \ref{main} we have that, for any $A>0$
\begin{equation}\label{ap1}
\B_q(t,s)=K(A;t)\int_0^{\infty/A}
\frac{x^{t-1}}{(1+x)_q^{t+s}}d_qx\ .
\end{equation}
By definition of $K(x,t)$ and using the results of Lemma \ref{prelim}
we get, after simple algebraic manipulations
\begin{equation}\label{ap2}
\frac{1}{x^t}(1+x)_q^{t+s}
=K\Big(\frac{1}{x};t\Big)\Big(1+\frac{q}{q^tx}\Big)_q^t(1+q^tx)_q^s\ .
\end{equation}
Since by Lemma \ref{prop_k} we have
$$
K\Big(\frac{1}{x};t\Big)=K(A;t)\ ,\quad \forall\ x=\frac{q^n}{A}\ ,\ n\in\Z\ ,
$$
when we substitute (\ref{ap2}) back into (\ref{ap1}) we get,
after a change of variable $y=q^tx$,
\begin{equation}\label{symmetric}
\B_q(t,s)=\int_0^{\infty/\alpha}
\frac{1}{y\Big(1+\frac{q}{y})_q^t(1+y)_q^s}d_qy\ ,\quad \forall\ \alpha>0\ .
\end{equation}
To conclude, we just notice that this integral expression of $\B_q(t,s)$
is manifestly symmetric in $t$ and $s$,
since performing the change of variable $x=\frac{q}{y}$ (namely 
applying the reciprocity relation (\ref{reciprocity})
gives the same integral with $\alpha$ replaced by $1/\alpha$
and $s$ replaced by $t$.


\ssection{Application 2: translation invariance of a certain type 
of improper integrals}

One of the main failures of the Jackson integral is 
that there is no analogue of the translation invariance identity
$$
\int_0^af(x)dx=\int_c^{a+c}f(x-c)dx\ ,
$$
obviously true for ``classical'' integrals.
By using Theorem \ref{main} we are able to write a $q$ analogue of
translation invariance for improper integrals of a special class of function,
namely of the form $x^\alpha/(1+x)_q^\beta$.
More precisely we want to prove the following
\begin{corollary}
For $\alpha>0$ and $\beta>\alpha+1$ we have
\begin{equation}\label{tr_inv}
\int_0^{\infty/\alpha}\frac{x^\alpha}{(1+x)_q^\beta}d_qx
=\frac{q}{q^\beta K(A,\alpha)}\int_1^{\infty/1}
\frac{x^\alpha\Big(1-\frac{1}{x}\Big)_q^\alpha}{x^\beta}d_qx\ .
\end{equation}
\end{corollary}
\begin{remark}
In the ``classical'' limit $q=1$, the right--hand side is obtained from the
left--hand side by translating $x\rightarrow x-1$.
\end{remark}
\begin{proof}
From the definition of $\B_q(t,s)$ we have
\begin{eqnarray}\label{cor1}
\B_q(t,s) &=& \int_0^1x^{s-1}(1-qx)_q^{t-1}d_qx \nonumber\\
 &=& \int_q^{\infty/1}
  \frac{1}{x^{s+1}}\Big(1-\frac{q}{x}\Big)_q^{t-1}d_qx \nonumber\\
 &=& \frac{1}{q^s}\int_1^{\infty/1}
  \frac{x^{t-1}\Big(1-\frac{1}{x}\Big)_q^{t-1}}{x^{t+s}}d_qx\ .
\end{eqnarray}
The first identity was obtained by applying (\ref{reciprocity})
and the second by a change of variable $y=x/q$.
From Theorem \ref{main} we also have
\begin{equation}\label{cor2}
\B_q(t,s)=K(A;t)\int_0^{\infty/A}
\frac{x^{t-1}}{(1+x)_q^{t+s}}d_qx\ .
\end{equation}
Equation (\ref{tr_inv}) is obtained by comparing (\ref{cor1}) and (\ref{cor2}),
after letting $\alpha=t-1,\ \beta=t+s$ and using the fact that
$K(A;\alpha+1)=q^\alpha K(A;\alpha)$.
\end{proof}


\ssection{Application 3: identities}

If we rewrite equations (\ref{t1}), (\ref{t2}) and (\ref{symmetric})
using the definition of improper integrals, we get some interesting identities
involving $q$--bilateral series.

After using the infinite product expansion (\ref{intr7}) 
of the $q$--exponential  function $e_q^x$, the expression (\ref{a}) 
of the $q$--gamma function,
the definition (\ref{intr15}) of the improper Jackson integral
and simple algebraic manipulations,
we can rewrite equation (\ref{t1}) as
\begin{equation}\label{jac1}
(1-q)_q^\infty(1+q^t/A)_q^\infty(1+qA/q^t)_q^\infty
=(1+qA)_q^\infty(1-q^t)_q^\infty
\sum_{n=-\infty}^\infty q^{tn}(1+1/A)_q^n\ .
\end{equation}
If we let $x=-q^t/A$ in equation (\ref{jac1}) we get
\begin{equation}\label{jac2}
(1-q)_q^\infty(1-x)_q^\infty(1-q/x)_q^\infty
=(1+qA)_q^\infty(1+Ax)_q^\infty
\sum_{n=-\infty}^\infty (-x)^n A^n(1+1/A)_q^n\ .
\end{equation}
This is a 1--parameter family of identities for the Jacobi triple product
$(1-q)_q^\infty(1-x)_q^\infty(1-q/x)_q^\infty\ ,$
parametrized by $A$. Notice that 
$$
\lim_{A\rightarrow0}A^n(1+1/A)_q^n=q^{n(n-1)/2}\ .
$$
This implies that, in the special case $A=0$, equation (\ref{jac2})
reduces to the famous Jacobi triple product identity (\ref{intr17a}).

Let's consider now equation (\ref{t2}).
After using the definition (\ref{intr15}) 
of improper $q$--integral,
the expression (\ref{intr9}) for the $q$--beta function
and simple algebraic manipulations, we can rewrite it as
\begin{equation}\label{ram1}
\sum_{n=-\infty}^\infty
\frac{(1+1/A)_q^n}{(1+q^{t+s}/A)_q^n}q^{tn}
=\frac{
(1-q)_q^\infty(1-q^{t+s})_q^\infty(1+q^t/A)_q^\infty(1+qA/q^t)_q^\infty
}{
(1+q^{t+s}/A)_q^\infty(1+qA)_q^\infty(1-q^t)_q^\infty(1-q^s)_q^\infty
}\ .
\end{equation}
Notice that, after letting $a=-1/A,\ b=-q^{t+s}/A,\ x=q^t$,
equation (\ref{ram1}) is equivalent to the famous 
Ramanujan's identity (\ref{intr17b}).
In other words, the proof of Theorem \ref{main} in Section 4 can be viewed 
as a new, more conceptual proof of Ramanujan's identity.

Finally we can rewrite equation (\ref{symmetric}) as
\begin{eqnarray}\label{sym1}
&& \sum_{n=-\infty}^\infty
\frac{
(1+1/\alpha)_q^n(1+q\alpha)_q^{-n}
}{
(1+q^s/\alpha)_q^n(1+q^{t+1}\alpha)_q^{-n}
} \nonumber\\
&& \hspace{-1cm} =\ \frac{
(1+1/\alpha)_q^\infty(1+q\alpha)_q^\infty
}{
(1+q^s/\alpha)_q^\infty(1+q^{t+1}\alpha)_q^\infty
}
\frac{
(1-q)_q^\infty(1-q^{t+s})_q^\infty
}{
(1-q^s)_q^\infty(1-q^t)_q^\infty}\ .
\end{eqnarray}
After letting $a=-1/\alpha,\ b=-q^s/\alpha,\ c=-q^{t+1}\alpha$,
equation (\ref{sym1}) reduces to (\ref{intr17c}).

\end{document}